\documentclass[11pt,a4paper,twoside]{article}

\usepackage{iftex}
\usepackage[T1]{fontenc}
\ifPDFTeX
  \usepackage[utf8]{inputenc}
\fi
\usepackage{lmodern}

\usepackage{amsmath,amsthm,amssymb,mathtools}
\usepackage{esint}

\usepackage{geometry}
\usepackage{xcolor}
\usepackage{microtype}
\usepackage{fancyhdr}
\usepackage{cite}
\usepackage[
  colorlinks=true,
  linkcolor=blue!45!black,
  citecolor=blue!45!black,
  urlcolor=blue!45!black
]{hyperref}

\allowdisplaybreaks[3]
\newtheoremstyle{paperplain}{0.65em}{0.65em}{\itshape}{}%
  {\bfseries}{.}{0.5em}{}
\theoremstyle{paperplain}
\newtheorem{theorem}{Theorem}[section]

\newtheorem{lemma}[theorem]{Lemma}

\newtheorem{definition}[theorem]{Definition}
\newtheoremstyle{paperremark}{0.65em}{0.65em}{\normalfont}{}%
  {\bfseries}{.}{0.5em}{}
\theoremstyle{paperremark}

\newcommand{\runningauthors}{Guohuan Qiu, Jin Yan, Xingchen Zhou, and Ruixuan Zhu}
\newcommand{\runningtitle}{}
\let\originaltitle\title
\renewcommand{\title}[2][]{%
  \originaltitle{#2}%
  \if\relax\detokenize{#1}\relax
    \gdef\runningtitle{#2}%
  \else
    \gdef\runningtitle{#1}%
  \fi
}

\numberwithin{equation}{section}
\newcommand{\R}{\mathbb R}
\DeclareMathOperator{\Sym}{Sym}

\DeclareMathOperator{\dom}{dom}

\DeclareMathOperator{\tr}{tr}
\DeclareMathOperator{\dist}{dist}
\title[Rigidity for Hessian quotients]{\textbf{Rigidity for semiconvex entire solutions\protect\\ of Hessian quotient equations}}

\author{\textbf{Guohuan Qiu}\and\textbf{Jin Yan}\and\textbf{Xingchen Zhou}\and\textbf{Ruixuan Zhu}}
\date{}

\makeatletter
\newcommand{\printfirstpagebottommatter}{%
  \begingroup
  \renewcommand{\thefootnote}{}%
  \renewcommand{\@makefntext}[1]{\noindent##1}%
  \footnotetext[0]{%
    \textit{Keywords.}
    Hessian quotient equation; Liouville theorem;
    semiconvex entire solution; Legendre transform.
    \par\smallskip
    \textit{2020 Mathematics Subject Classification.}
    Primary 35B53; Secondary 35J60.%
  }%
  \endgroup
}
\makeatother

\newcommand{\printauthorinformation}{%
  \par\bigskip
  \begingroup
  \small
  \setlength{\parindent}{0pt}%

  \noindent
  \textbf{Guohuan Qiu}\par
  State Key Laboratory of Mathematical Sciences,
  Academy of Mathematics and Systems Science,
  Chinese Academy of Sciences,
  Beijing 100190, China; \\
  Institute of Mathematics,
  Academy of Mathematics and Systems Science,
  Chinese Academy of Sciences,
  55 Zhongguancun East Road, Beijing 100190, China.\par
  \textit{Email:}
  \href{mailto:qiugh@amss.ac.cn}{qiugh@amss.ac.cn}.

  \par\medskip

  \noindent
  \textbf{Jin Yan}\par
  Institute of Mathematics,
  Academy of Mathematics and Systems Science,
  Chinese Academy of Sciences,
  55 Zhongguancun East Road, Beijing 100190, China.\par
  \textit{Email:}
  \href{mailto:yanjin@amss.ac.cn}{yanjin@amss.ac.cn}.

  \par\medskip

  \noindent
  \textbf{Xingchen Zhou}\par
  School of Mathematics and Statistics,
  Hainan University,
  Haikou 570228, China.\par
  \textit{Email:}
  \href{mailto:zxc3zxc4zxc5@stu.xjtu.edu.cn}{zxc3zxc4zxc5@stu.xjtu.edu.cn}.

  \par\medskip

  \noindent
  \textbf{Ruixuan Zhu}\par
  Institute for Theoretical Sciences,
  Westlake University,
  Hangzhou 310030, China.\par
  \textit{Email:}
  \href{mailto:zhuruixuan@westlake.edu.cn}{zhuruixuan@westlake.edu.cn}.

  \endgroup
}

\hypersetup{%
  pdftitle={Rigidity for semiconvex entire solutions of Hessian quotient equations},
  pdfauthor={Guohuan Qiu, Jin Yan, Xingchen Zhou, and Ruixuan Zhu},
  pdfsubject={Rigidity of semiconvex entire solutions of Hessian quotient equations},
  pdfkeywords={Hessian quotient equation, Liouville theorem, semiconvex entire solution, Legendre transform}
}

\AtBeginDocument{\nocite{BCD,CabreCaffarelli2003,CC,CNS,Chen,Ding,DZ,Guan,Gutierrez,HL,HLAE,QY,Rockafellar,Savin,SY,SYrigidity}}
\begin{document}

\maketitle
\printfirstpagebottommatter
\begin{abstract}
We prove that every smooth admissible semiconvex entire solution of $\sigma_k(D^2u)/\sigma_l(D^2u)=1$ in $\mathbb R^n$,
where $0\le l<k\le n$ and $k-l=1$ or $2$,
is a quadratic polynomial.
The proof uses a shifted Legendre transform and interior Hessian
estimates for Hessian quotient equations to reduce the problem
to equations of lower order. We then conclude by induction on the order of the equation.
\end{abstract}

\section{Introduction}

Let $\sigma_j$ denote the elementary symmetric polynomial
of degree $j$, with $\sigma_0=1$, and let
\[
 \Gamma_k^n=\{\lambda\in\R^n:\sigma_j(\lambda)>0,
 \ 1\le j\le k\}
\]
be the G{\aa}rding cone. For a symmetric matrix, $\sigma_j$ is evaluated
at its eigenvalues. A smooth function is called $k$-admissible if
$\lambda(D^2u)\in\Gamma_k^n$; we also write $D^2u\in\Gamma_k^n$.
We consider entire admissible solutions of the Hessian quotient equation
\[
 \frac{\sigma_k(D^2u)}{\sigma_l(D^2u)}=1
 \qquad\text{in }\R^n,\qquad 0\le l<k\le n.
\]
The case $l=0$ is the $k$-Hessian equation. The operator
$(\sigma_k/\sigma_l)^{1/(k-l)}$ is elliptic and concave on
$\Gamma_k^n$; see \cite{CNS,Spruck}.

The rigidity problem asks under what conditions an entire solution
must be a quadratic polynomial. In the linear case \(k=1\),
semiconvexity bounds the Hessian, and rigidity follows from the
Liouville theorem for harmonic functions.
For the Monge--Amp\`ere equation, the theorem of
J{\"o}rgens \cite{Jorgens}, Calabi \cite{Calabi}, and
Pogorelov \cite{Pogorelov} states that every convex entire solution
is quadratic.

Related rigidity results concern the special Lagrangian equation
\begin{equation*}
\sum_{i=1}^n\arctan\lambda_i(D^2u)=\Theta.
\end{equation*}
For smooth entire solutions, Yuan \cite{Yuan2002} proved rigidity
under convexity in every dimension and under semiconvexity
in dimension three. He also established rigidity in the
supercritical phase range \(|\Theta|>(n-2)\pi/2\)
\cite{Yuan2006}. Li--Li--Yuan \cite{LiLiYuan2020}
extended the semiconvex rigidity result to dimension four
for arbitrary phase. Ogden--Yuan \cite{OgdenYuan2025}
obtained further rigidity results in higher dimensions
under bounds on the smallest or largest Hessian eigenvalue
depending on the phase.
Mooney--Shankar \cite{MooneyShankar2025} also proved rigidity
for entire viscosity solutions with
\(-(n-2)\pi/2<\Theta<\pi/2\), under an explicit
semiconvexity bound depending on the phase.

For the \(2\)-Hessian equation, Chang--Yuan \cite{CY} proved
that every smooth entire solution of \(\sigma_2(D^2u)=1\)
is quadratic under the almost convexity condition
\(D^2u\ge\left(\delta-\sqrt{\frac{2}{n(n-1)}}\right)I\)
for some \(\delta>0\). Shankar--Yuan \cite{SYrigidity} later extended this rigidity
result to all smooth semiconvex entire solutions, settling the \(\sigma_2\) case of the conjecture of Chang--Yuan.
More recently, Zhou--Zhu \cite{ZZ} proved that every convex
\(C^2\) entire solution of \(\sigma_k(D^2u)=1\) is quadratic
for \(3\le k\le n-1\).

There are also many results under a quadratic lower bound condition \(u(x)\ge a|x|^2-b\) with \(a>0\). Earlier rigidity results required additional convexity conditions or restrictions on \(k\) \cite{BCGJ,LRW,ChuDinew,Tu2024,Zhang2025}.
Recently, Wang \cite{Wang2026} treated \(k\)-admissible solutions
with \(n/2<k<n\), and Hu--Jiao \cite{HuJiao2026} proved rigidity for all \(2\le k\le n\) under the growth condition alone.

For Hessian quotient equations,
Bao--Chen--Guan--Ji \cite{BCGJ} proved rigidity for strictly
convex entire solutions of \(\sigma_n/\sigma_l=1\)
under a quadratic upper growth bound.
They also observed that the growth assumption can be removed
when \(l=n-1\).
Using the rigidity theorem of Shankar--Yuan,
Lu--Sroka \cite{LuSroka} proved rigidity for
\(\sigma_2/\sigma_1=1\) under semiconvexity alone
by subtracting a fixed quadratic function and reducing
the equation to the \(\sigma_2\) equation.
See also \cite{Du2022,MeiYan2026,LiWu2026} for related
rigidity results.

However, a positive quadratic lower bound excludes
admissible quadratic solutions with indefinite Hessians.
Semiconvexity allows such solutions and is therefore a natural
condition for rigidity.
On the other hand, \(k\)-admissibility alone does not imply
rigidity in general.
For example, Warren \cite{WarrenCounterexamples} constructed smooth
nonpolynomial entire \(k\)-admissible solutions of
\(\sigma_k(D^2v)=1\) in \(\R^n\) for every \(k\ge2\)
and \(n\ge2k-1\).
For \(k=2\) and \(n=3\), one such solution is
\begin{equation*}
 v(x_1,x_2,t)
 = (x_1^2+x_2^2-1)e^t+\frac14e^{-t}.
\end{equation*}

In this paper, we prove that semiconvexity alone implies
rigidity for admissible entire solutions of Hessian quotient
equations with \(k-l=1\) or \(2\), without any growth
assumption at infinity.

\begin{theorem}\label{thm:main}
Let $0\le l<k\le n$ and $u\in C^\infty(\R^n)$ satisfy
\begin{equation}\label{eq:main-equation}
 D^2u\in\Gamma_k^n,\qquad
 \frac{\sigma_k(D^2u)}{\sigma_l(D^2u)}=1,\qquad D^2u\ge-KI
\end{equation}
for some finite $K\ge0$. If $k-l=1$ or $2$, then $u$ is a quadratic polynomial.
\end{theorem}

Our proof builds on the work of Shankar--Yuan \cite{SYrigidity}.
Their Legendre--Lewy transform produces a potential with bounded
Hessian satisfying a uniformly elliptic equation that need not
be concave.
They establish a global Jacobi inequality for
\(\log(\Delta u+J)\), with a suitable constant \(J\), and show
that \((\Delta u+J)^{-1/3}\) remains a positive supersolution
for the transformed linearized operator.
The weak Harnack inequality then shows that this quantity
is nearly constant in measure on a suitable ball.
Through the transformed equation, they transfer this information
to a concave, uniformly elliptic Hessian quotient.
Comparison with a solution of the corresponding Dirichlet
problem, together with the ABP and Evans--Krylov estimates,
yields approximation by a quadratic polynomial.
They then adapt the iteration method of Caffarelli--Yuan
\cite{CaffarelliYuan2000} to improve this approximation at
smaller scales and obtain a uniform interior H\"older estimate
for the transformed Hessian.
A scaling argument finally shows that both the transformed
and original Hessians are constant.

For the quotient equations considered here, we combine the
Legendre transform and the weak Harnack inequality with a
reduction to entire solutions of lower-order equations.
We first treat adjacent quotients and write
\begin{equation*}
 F_{n,k}=\frac{\sigma_k}{\sigma_{k-1}}.
\end{equation*}
We proceed by induction on \(k\), simultaneously in all
dimensions \(n\ge k\).
If \(D^2u\) is bounded, rigidity follows from the
Evans--Krylov estimate and scaling.
It therefore remains to exclude an unbounded Hessian.
The argument has three steps.

\medskip
\noindent\textit{1. A limit with vanishing Hessian determinant.}
Let \(\kappa=K+1\) and define
\begin{equation*}
 w=\left(u+\frac{\kappa}{2}|x|^2\right)^*.
\end{equation*}
Then \(w\) is smooth on \(\R^n\) and satisfies
\(0<D^2w\le I\).
Following Shankar--Yuan \cite{SY,SYrigidity}, we combine
a Jacobi inequality with the Legendre transform and
the weak Harnack inequality.
The required ellipticity and quantitative concavity estimates
are due to Dong--Zhang \cite{DZ}.
Inspired by Ding \cite{Ding}, we use the assumed unboundedness
of \(D^2u\) to show that the averages of \(\det D^2w\)
over expanding balls tend to zero.
Normalized quadratic rescalings \(w_j\) then converge
locally in \(C^1\) to a convex function
\(W\in C^{1,1}_{\mathrm{loc}}(\R^n)\) satisfying
\begin{equation*}
 0\le D^2W\le I,
 \qquad \det D^2W=0
 \quad\text{almost everywhere}.
\end{equation*}

\medskip
\noindent\textit{2. Recovering an entire solution of lower order.}
Since \(W\) has zero Monge--Amp\`ere measure, convex duality
shows that it is constant along fixed directions.
Removing all these directions yields a splitting
\begin{equation*}
 W(y'',y')=V(y'),
 \qquad (y'',y')\in\R^m\times\R^{n-m},
\end{equation*}
where \(V\) has no nonzero global affine direction.
Our description of the possible limiting Hessians shows
that \(1\le m\le k-1\) and identifies a uniformly elliptic
equation satisfied by \(V\).

The main step is to recover a smooth entire solution from \(V\).
The splitting alone does not ensure that \(V^*\) is finite
everywhere, and \(D^2V\) may still have zero eigenvalues.
We first use the uniformly elliptic equation for \(V\)
to prove that \(V^*\) is finite on all of \(\R^{n-m}\).
If its domain were a proper subset, the Harnack inequality
for translation differences of \(V\) would force a nonzero
global affine direction, contradicting the choice of \(V\).
We then regularize \(V\) before taking its conjugate and
pass to the limit to show that
\begin{equation*}
 z(x)=V^*(x)-\frac{\kappa}{2}|x|^2
\end{equation*}
is an entire admissible viscosity solution of
\begin{equation*}
 \frac{\sigma_{k-m}(D^2z)}
      {\sigma_{k-m-1}(D^2z)}=1
 \quad\text{in }\R^{n-m},
 \qquad D^2z\ge-KI.
\end{equation*}
The interior Hessian estimates of Qiu--Yan
\cite[Theorem~1.1]{QY}, applied to smooth Dirichlet
approximations, imply that \(z\) is smooth.
Thus we recover an entire solution of an adjacent quotient
equation of lower order while preserving the original
semiconvexity bound.

\medskip
\noindent\textit{3. Induction and the final contradiction.}
Since \(1\le k-m<k\), the induction hypothesis implies that
\(z\) is quadratic, and hence so are \(V\) and \(W\).
Savin's small perturbation theorem \cite{Savin} then
upgrades the convergence of \(w_j\) to convergence of
their Hessians near the origin.
Quadratic rescaling leaves \(D^2w(0)\) unchanged,
so \(D^2w(0)=D^2W\).
This contradicts \(D^2w(0)>0\) and \(\det D^2W=0\),
and rules out an unbounded Hessian.

\medskip
The case \(k-l=2\) follows by applying the adjacent quotient
result to \(U(x,t)=u(x)+t^2/2\), as in \cite{QY}.

The paper is organized as follows. Section~\ref{sec:prelim} collects several results needed for the proof. Section~\ref{sec:spectral} studies the transformed equation and the
closure of its admissible Hessians. Section~\ref{sec:jacobi}
constructs the limit $W$. Section~\ref{sec:splitting} splits
this limit and uses convex conjugation to recover a smooth
solution of a lower-order equation. Section~\ref{sec:main-proof}
completes the induction and the
proof of Theorem~\ref{thm:main}.

\section{Preliminaries}\label{sec:prelim}
Throughout this paper, the semiconvexity constant
$K\ge0$ is fixed, and we set
\begin{equation}\label{eq:kappa}
 \kappa=K+1.
\end{equation}
Unless otherwise stated, $c$ and $C$ denote positive constants depending
only on the dimension, the order, and $K$. For $F=F_{n,k}$,
we write 
$$
F_i=\frac{\partial F}{\partial\lambda_i}\qquad \text{and} \qquad F_{ij}=\frac{\partial^2F}{\partial\lambda_i\partial\lambda_j}.
$$
We first recall the definition of an admissible viscosity solution; see \cite{CIL,Trudinger1990} for details.
\begin{definition}
A continuous function $z$ is an admissible
viscosity solution of $F_{n,k}(D^2z)=1$ if every upper test
$\varphi$ at $x_0$ satisfies
\[
 D^2\varphi(x_0)\in\Gamma_k^n,
 \qquad F_{n,k}(D^2\varphi(x_0))\ge1,
\]
and every lower test $\varphi$ at $x_0$ with
$D^2\varphi(x_0)\in\Gamma_k^n$ satisfies
\[
 F_{n,k}(D^2\varphi(x_0))\le1.
\]
\end{definition}
Next, we gather several known facts, which will be used repeatedly.
\begin{lemma}\label{lem:cone-newton}
Let \(1\le k\le n\). The G{\aa}rding cone \(\Gamma_k^n\) is open,
convex, and preserved by positive semidefinite additions.
The quotient \(F_{n,k}=\sigma_k/\sigma_{k-1}\) is smooth,
strictly elliptic, concave, and homogeneous of degree one
on \(\Gamma_k^n\).

If \(\eta\in\Gamma_k^n\) satisfies \(F_{n,k}(\eta)\ge c_0>0\), then
\begin{equation}\label{eq:newton-lower}
  \sigma_j(\eta)\ge\binom nj
  \left(\frac{kc_0}{n-k+1}\right)^j,
  \qquad 1\le j\le k.
\end{equation}
In particular,
\begin{equation}\label{eq:boundary-zero}
  F_{n,k}(\eta)\le\frac{n-k+1}{k}
  \left(\frac{\sigma_k(\eta)}{\binom nk}\right)^{1/k}.
\end{equation}
\end{lemma}

\begin{proof}
These follow directly from the generalized Newton--Maclaurin inequalities; see \cite[Lemma~2.10]{Spruck}. 
\end{proof}

\begin{lemma}\label{lem:DZ-estimates}
Let $2\leq k\leq n$, and suppose that
\[
 \lambda\in\Gamma_k^n,\qquad
 \lambda_1\ge\cdots\ge\lambda_n\ge-K,
 \qquad F_{n,k}(\lambda)=1.
\]
Then there are constants $c,C>0$, depending only on $n,k,K$, such that
\begin{equation}\label{eq:dual-ellipticity}
 c\le F_i(\lambda)(\lambda_i+\kappa)^2\le C,
 \qquad 1\le i\le n.
\end{equation}
There also exist $L_*>1$ and $\delta,\gamma>0$, with the same
dependence, such that, if $L=\lambda_1\ge L_*$, then
\begin{equation}\label{eq:tangent}
 -\sum_{i,j}F_{ij}\zeta_i\zeta_j
 +\frac{2}{(1+\delta)L}\sum_{i>1}F_i\zeta_i^2
 \ge(1+\gamma)\frac{F_1\zeta_1^2}{L}
\end{equation}
for every $\zeta\in\R^n$ satisfying $\sum_iF_i\zeta_i=0$.
\end{lemma}

\begin{proof}
These estimates follow from \cite[Lemmas~2.3 and~2.5]{DZ}
with \(F_{n,k}=1\).
\end{proof}

\begin{lemma}\label{lem:strong-viscosity}
Let \(v\in C^{1,1}_{\mathrm{loc}}(\Omega)\), where
\(\Omega\subset\R^d\) is open.
Suppose that \(\mathcal E\subset\Sym(d)\) is closed and satisfies
\(A+P\in\mathcal E\) whenever \(A\in\mathcal E\) and \(P\ge0\).
If \(D^2v\in\mathcal E\) almost everywhere in \(\Omega\), then
every upper test \(\varphi\) for \(v\) at \(x_0\in\Omega\)
satisfies \(D^2\varphi(x_0)\in\mathcal E\).
\end{lemma}

\begin{proof}
Since \(v\) is locally semiconvex, the assertion follows
from the almost-everywhere theorem of Harvey and Lawson
\cite[Theorem~4.1]{HLAE}.
\end{proof}

For $0<a\le b$, the Pucci extremal operators are defined by
\[
 \mathcal P^+_{a,b}(M)=b\tr M^+-a\tr M^- ,
 \qquad
 \mathcal P^-_{a,b}(M)=a\tr M^+-b\tr M^- ,
\]
where $M^+$ and $M^-$ are the positive and negative spectral parts
of $M$, so that $M=M^+-M^-$ and $M^\pm\ge0$.

\begin{lemma}\label{lem:harnack}
Let $q\in C(B_{2R})$ be nonnegative. If
\[
 \mathcal P^-_{a,b}(D^2q)\le0
 \qquad\text{in }B_{2R}
\]
in the viscosity sense, then there exist $s>0$ and $C<\infty$,
depending only on the dimension and $a,b$, such that
\begin{equation}\label{eq:prelim-weak-harnack}
 \left(\fint_{B_R}q^s\right)^{1/s}
 \le C\inf_{B_R}q.
\end{equation}
If, in addition,
$\mathcal P^+_{a,b}(D^2q)\ge0$ in $B_{2R}$, then
\begin{equation}\label{eq:prelim-harnack}
 \sup_{B_R}q\le C\inf_{B_R}q.
\end{equation}
\end{lemma}

\begin{proof}
Both inequalities are standard consequences of the Krylov--Safonov theory; see \cite[Chapter~4]{CC}.
\end{proof}

\begin{lemma}\label{lem:concave-regularity}
Let $H:\Sym(n)\to\R$ be continuous, concave, and uniformly
elliptic, with
\[
 a\tr B\le H(A+B)-H(A)\le b\tr B
 \qquad(A\in\Sym(n),\ B\ge0).
\]
If a bounded continuous viscosity solution $v$ satisfies
$H(D^2v)=c$ in $B_1$, then
\begin{equation}\label{eq:prelim-evans-krylov}
 \|v\|_{C^{2,\alpha}(B_{1/2})}
 \le C\bigl(\|v\|_{L^\infty(B_1)}+|c-H(0)|\bigr),
\end{equation}
where $\alpha\in(0,1)$ and $C$ depend only on $n,a,b$.
\end{lemma}

\begin{proof}
This is the classical Evans--Krylov estimate; see \cite[Chapter~6]{CC}.
\end{proof}

\begin{lemma}\label{lem:savin-small}
Let $G:\Sym(n)\to\R$ be continuous and globally uniformly
elliptic. Suppose that $G(0)=0$ and that $G$ is $C^2$ in a
neighborhood of zero, with bounded first and second derivatives
on a smaller neighborhood. There exist
$\varepsilon_*>0$, $\alpha\in(0,1)$, and $C<\infty$,
depending only on the dimension, the ellipticity constants,
and these local bounds and neighborhoods, such that every
continuous viscosity solution of
\[
 G(D^2v)=0\quad\text{in }B_1,
 \qquad \|v\|_{L^\infty(B_1)}\le\varepsilon_*
\]
satisfies
\begin{equation}\label{eq:prelim-savin}
 \|v\|_{C^{2,\alpha}(B_{1/2})}\le C.
\end{equation}
\end{lemma}

\begin{proof}
This is Savin's small perturbation theorem; see \cite[Theorem~1.3]{Savin}.
\end{proof}

\section{The transformed equation}\label{sec:spectral}
In this section, we use a shifted Legendre transform to obtain
a convex function with uniformly bounded Hessian.
We then describe the Hessians that can arise in the limit
and extend the transformed operator while keeping it
uniformly elliptic. These results will be used later to pass to limits and derive lower-order equations.

Let \(u\) be a smooth entire admissible solution of
\(F_{n,k}(D^2u)=1\) with \(D^2u\ge-KI\). Set
\begin{equation}\label{eq:fixed-legendre}
 v(x)=u(x)+\frac\kappa2|x|^2,
 \qquad y=Dv(x),\qquad w=v^*.
\end{equation}
Here the Legendre transform is defined by
$$
w(y)=v^*(y)=\sup_{x\in\R^n}\{x\cdot y-v(x)\}.
$$
Since \(D^2v\ge I\), the conjugate \(w\) is well-defined on \(\R^n\), and satisfies the following properties:
\begin{equation}\label{eq:dual-hessian}
 Dw(y)=x,\qquad
 D^2w(y)=(D^2u(x)+\kappa I)^{-1},\qquad
 0<D^2w\le I.
\end{equation}
In particular, we have
$$
0\le w(y)-w(0)-Dw(0)\cdot y\le\frac12|y|^2.
$$
For any dimension \(d\) and \(1\le p\le d\), define
\begin{equation}\label{eq:physical-set}
 \Sigma_{d,p}
 =\overline{\left\{(A+\kappa I)^{-1}:
 A\in\Gamma_p^d,\quad F_{d,p}(A)=1,\quad A\ge-KI\right\}}.
\end{equation}
It is easy to see that the sets \(\Sigma_{d,p}\) are compact, and every
\(C\in\Sigma_{d,p}\) satisfies \(0\le C\le I\).
The next lemma describes all matrices in these sets, including
singular ones.

\begin{lemma}\label{lem:strata}
A matrix \(C\) belongs to \(\Sigma_{n,k}\) if and only if,
after an orthogonal change of coordinates,
$$
C=0_r\oplus B,\qquad B>0,\qquad 0\le r\le k-1,
$$
and the spectrum \(\eta\) of \(B^{-1}-\kappa I\) satisfies
\begin{equation}\label{eq:tail-equation}
 \eta\in\Gamma_{k-r}^{n-r},\qquad
 F_{n-r,k-r}(\eta)=1,\qquad \eta_i\ge-K.
\end{equation}
\end{lemma}

\begin{proof}
We first show that the \(k\)th largest eigenvalue is bounded.
For large \(L\), consider
$$
\beta_L=(\underbrace{L,\ldots,L}_{k},
                \underbrace{-K,\ldots,-K}_{n-k}).
$$
If an ordered admissible spectrum \(\lambda\) satisfies
\(\lambda_i\ge-K\), \(F_{n,k}(\lambda)=1\), and
\(\lambda_k\ge L\), then \(\lambda_i\ge(\beta_L)_i\)
for every \(i\). By monotonicity,
\(1=F_{n,k}(\lambda)\ge F_{n,k}(\beta_L)>1\), a contradiction.

Let \(C_j\to C\), where each \(C_j\) belongs to the set in
\eqref{eq:physical-set} before taking the closure, and suppose
\(\dim\ker C=r\). Then exactly \(r\) eigenvalues of
\(C_j^{-1}-\kappa I\) tend to \(+\infty\).
The remaining eigenvalues converge to the spectrum \(\eta\)
of \(B^{-1}-\kappa I\), where \(B\) is the restriction of
\(C\) to its range. The bound above shows that \(r\le k-1\).
Set
$$
d=n-r,\qquad p=k-r,
$$
and let \(\eta^{(j)}\) be the vector of the remaining eigenvalues.
The Newton--Maclaurin inequalities imply
\(\eta^{(j)}\in\Gamma_p^d\) and
\(F_{d,p}(\eta^{(j)})\ge1\). Passing to the limit, we obtain \(\eta\in\Gamma_p^d\).

Let \(a^{(j)}=(a_1^{(j)},\ldots,a_r^{(j)})\) denote the
eigenvalues that tend to infinity, so
\(a_i^{(j)}\to+\infty\) for every \(i\).
For \(s=k-1,k\) and sufficiently large \(j\),
\begin{align*}
 \frac{\sigma_s(a^{(j)},\eta^{(j)})}
      {\prod_{i=1}^r a_i^{(j)}}
 &=\sum_{q=0}^r
 \sigma_q\!\left(\frac1{a_1^{(j)}},\ldots,
                 \frac1{a_r^{(j)}}\right)
 \sigma_{s-r+q}(\eta^{(j)})\\
 &=\sigma_{s-r}(\eta^{(j)})+o(1).
\end{align*}
Dividing
\(\sigma_k(a^{(j)},\eta^{(j)})
 =\sigma_{k-1}(a^{(j)},\eta^{(j)})\)
by \(\prod_{i=1}^r a_i^{(j)}\) and passing to the limit, we obtain
$$
\sigma_p(\eta)=\sigma_{p-1}(\eta).
$$
Since \(\eta\in\Gamma_p^d\), we have \(F_{d,p}(\eta)=1\) with \(\eta_i\ge-K\).

Conversely, suppose that \(C=0_r\oplus B\) satisfies the stated
conditions. For \(r=0\), the claim follows from the definition.
For \(r\ge1\), set
$$
\xi(T,t)=(\underbrace{T,\ldots,T}_{r},\eta+t\mathbf1),
\qquad t\ge0.
$$
When \(T\) is sufficiently large, \(\xi(T,0)\) lies in
\(\Gamma_k^n\), and increasing \(t\) keeps the vector in this cone.
For each fixed \(t\ge0\),
$$
F_{n,k}(\xi(T,t))
=\frac{\sigma_p(\eta+t\mathbf1)+O(T^{-1})}
      {\sigma_{p-1}(\eta+t\mathbf1)+O(T^{-1})}
\longrightarrow F_{d,p}(\eta+t\mathbf1).
$$
At \(t=0\), the function \(T\mapsto F_{n,k}(\xi(T,0))\)
is strictly increasing by ellipticity. Its limit is
\(F_{d,p}(\eta)=1\), so
$$
F_{n,k}(\xi(T,0))<1
$$
when \(T\) is sufficiently large. Strict ellipticity also implies
$$
F_{d,p}(\eta+\mathbf1)>F_{d,p}(\eta)=1.
$$
We may therefore choose \(T_1>0\) so that, whenever \(T\ge T_1\),
we have \(\xi(T,t)\in\Gamma_k^n\) for all \(t\ge0\) and
$$
F_{n,k}(\xi(T,0))<1<F_{n,k}(\xi(T,1)).
$$
By continuity and strict monotonicity in \(t\), there is a unique
\(t_T\in(0,1)\) such that
$$
F_{n,k}(\xi(T,t_T))=1.
$$
Fix \(\delta\in(0,1)\). Since
$$
F_{n,k}(\xi(T,\delta))
\longrightarrow F_{d,p}(\eta+\delta\mathbf1)>1,
$$
there exists \(T_\delta\ge T_1\) such that
$$
F_{n,k}(\xi(T,\delta))
>F_{n,k}(\xi(T,t_T))=1
\qquad\text{for all }T\ge T_\delta.
$$
Strict monotonicity then implies \(0<t_T<\delta\) whenever
\(T\ge T_\delta\). Since \(\delta\in(0,1)\) was arbitrary,
we conclude that \(t_T\to0\) as \(T\to\infty\).

Define
$$
A_T=TI_r\oplus\bigl(B^{-1}-\kappa I_d+t_T I_d\bigr).
$$
Its eigenvalues are \(\xi(T,t_T)\), so
\(A_T\in\Gamma_k^n\) and \(F_{n,k}(A_T)=1\).
Since \(\eta_i\ge-K\) and \(t_T>0\), we also have
\(A_T\ge-KI_n\), and
$$
(A_T+\kappa I_n)^{-1}
=\frac{I_r}{T+\kappa}
 \oplus(B^{-1}+t_T I_d)^{-1}
\longrightarrow 0_r\oplus B=C.
$$
Thus \(C\in\Sigma_{n,k}\).
\end{proof}

For a positive definite dual Hessian \(C\), the transformed
equation is
$$
F_{n,k}(C^{-1}-\kappa I)=1.
$$
To handle singular limits, we write the equation using polynomials
in \(C\). Define \(T_j:\Sym(n)\to\R\) by
\begin{equation}\label{eq:generating-polynomial}
 \det\bigl(C+t(I-\kappa C)\bigr)
 =\sum_{j=0}^nT_j(C)t^j.
\end{equation}

For \(C>0\),
$$
C+t(I-\kappa C)=C\bigl(I+t(C^{-1}-\kappa I)\bigr),
$$
so
$$
\det\bigl(C+t(I-\kappa C)\bigr)
=\det C\sum_{j=0}^n\sigma_j(C^{-1}-\kappa I)t^j.
$$
Comparing the corresponding coefficients in the above equality, we obtain
$$
T_j(C)=\det C\,\sigma_j(C^{-1}-\kappa I).
$$
Hence, when \(C>0\) and \(C^{-1}-\kappa I\in\Gamma_k^n\), we have
$$
\frac{T_k(C)}{T_{k-1}(C)}=F_{n,k}(C^{-1}-\kappa I).
$$
Set
\begin{equation}\label{eq:dual-polynomial-operator}
 H_0(C)=1-\frac{T_k(C)}{T_{k-1}(C)}
\end{equation}
wherever \(T_{k-1}(C)\ne0\). Then for any positive definite \(C\in\Sigma_{n,k}\), the transformed equation is exactly \(H_0(C)=0\).

The next lemma shows that \(T_{k-1}>0\) on \(\Sigma_{n,k}\)
and extends \(H_0\) while keeping it uniformly elliptic.
We will use this extension to study the lower-dimensional limit
and to apply Savin's theorem at the end of the proof.

\begin{lemma}\label{lem:operator}
The denominator in \eqref{eq:dual-polynomial-operator} is strictly
positive on \(\Sigma_{n,k}\). The function \(H_0\) is smooth and
uniformly elliptic in a fixed neighborhood of \(\Sigma_{n,k}\),
and \(H_0=0\) on \(\Sigma_{n,k}\). Moreover,
$$
C,D\in\Sigma_{n,k},\quad C\le D\quad\Longrightarrow\quad C=D.
$$
There also exists a continuous operator \(\widehat H\) on \(\Sym(n)\)
that agrees with \(H_0\) in a neighborhood of \(\Sigma_{n,k}\)
and is uniformly elliptic on all of \(\Sym(n)\).
\end{lemma}

\begin{proof}
For \(C=0_r\oplus B\) as in Lemma~\ref{lem:strata},
$$
\det\bigl(C+t(I-\kappa C)\bigr)
=t^r\det\bigl(B+t(I-\kappa B)\bigr).
$$
Thus
$$
T_{k-1}(C)=\det B\,\sigma_{k-r-1}(B^{-1}-\kappa I)>0,
\qquad T_k(C)=T_{k-1}(C).
$$
Since \(\Sigma_{n,k}\) is compact, \(T_{k-1}\) has a positive
lower bound there. Hence \(H_0\) is smooth near \(\Sigma_{n,k}\)
and vanishes on \(\Sigma_{n,k}\).

First take \(C\) in the set in \eqref{eq:physical-set} before
the closure, and set \(A=C^{-1}-\kappa I\).
For any symmetric matrix \(N\),
$$
DH_0(C)[N]=\tr\bigl(C^{-1}DF_{n,k}(A)C^{-1}N\bigr).
$$
In a common eigenbasis of \(A\) and \(C\), the coefficient matrix
has eigenvalues
$$
F_i(\lambda)(\lambda_i+\kappa)^2.
$$
For \(2\le k<n\), \eqref{eq:dual-ellipticity} bounds these
eigenvalues uniformly above and below by positive constants.
For \(k=n\), the same bounds follow by direct computation.
By continuity, they also hold on \(\Sigma_{n,k}\).
Compactness then allows us to choose \(a_0,b_0,s>0\) such that
\begin{equation}\label{eq:local-operator-ellipticity}
 a_0I\le DH_0(X)\le b_0I
 \quad\text{when }\dist(X,\Sigma_{n,k})<s.
\end{equation}

To prove the order property, let \(C,D\in\Sigma_{n,k}\) satisfy
\(C\le D\). Then clearly \(\ker D\subset\ker C\).
Restrict both matrices to \(L=(\ker D)^\perp\), and denote
the restrictions by \(C_L,D_L\). Then \(D_L>0\).
Set \(r_D=\dim\ker D\) and \(p=k-r_D\).
By Lemma~\ref{lem:strata}, we know
$$
A_0=D_L^{-1}-\kappa I\in\Gamma_p^{n-r_D},
\qquad F_{n-r_D,p}(A_0)=1.
$$
Now suppose the contrary that \(C\ne D\). For \(0\le t<1\), set
$$
B_t=(1-t)D_L+tC_L,\qquad A_t=B_t^{-1}-\kappa I.
$$
Since
$$
A_t'=B_t^{-1}(D_L-C_L)B_t^{-1}\ge0,
$$
we know \(A_t\ge A_0\), so the path stays in
\(\Gamma_p^{n-r_D}\). Therefore, strict ellipticity implies
$$
\frac{d}{dt}F_{n-r_D,p}(A_t)
=\tr\bigl(DF_{n-r_D,p}(A_t)A_t'\bigr)>0.
$$
Let \(r_C=\dim\ker C\) and \(C_+=C|_{(\ker C)^\perp}\).
As \(t\uparrow1\), exactly \(r_C-r_D\) eigenvalues of \(A_t\)
tend to infinity, while the others converge to those of
\(C_+^{-1}-\kappa I\). The same expansion used in the proof
of Lemma~\ref{lem:strata} shows that
$$
\lim_{t\uparrow1}F_{n-r_D,p}(A_t)
=F_{n-r_C,k-r_C}(C_+^{-1}-\kappa I)=1.
$$
But strict monotonicity implies, for any fixed \(t_0\in(0,1)\),
$$
\lim_{t\uparrow1}F_{n-r_D,p}(A_t)
\ge F_{n-r_D,p}(A_{t_0})
>F_{n-r_D,p}(A_0)=1,
$$
a contradiction. Thus \(C=D\).

We now construct the global extension. Write a symmetric matrix
as \(E=E^+-E^-\), where \(E^+,E^-\) are its positive and negative
parts, and recall
$$
\mathcal P^+_{a,b}(E)=b\tr E^+-a\tr E^-,
\qquad
\mathcal P^-_{a,b}(E)=a\tr E^+-b\tr E^-.
$$
We use the infimum formula in
\cite[Remark~1.5]{CabreCaffarelli2003} to extend \(H_0\) from a
neighborhood of \(\Sigma_{n,k}\) to \(\Sym(n)\).
We first prove the bounds needed for this extension.

The order property shows that \(\tr(C-D)^+\) and \(\tr(C-D)^-\)
are both positive for distinct \(C,D\in\Sigma_{n,k}\).
Set \(a=a_0\). By compactness, we can choose \(b\ge b_0\)
sufficiently large that
\begin{equation}\label{eq:far-pair-margin}
 \mathcal P^+_{a,b}(C-D)\ge1,
 \qquad
 \mathcal P^-_{a,b}(C-D)\le-1
 \quad\text{whenever }C,D\in\Sigma_{n,k},\quad |C-D|\ge s/4.
\end{equation}
Choose \(\varepsilon>0\) such that
$$
\varepsilon<s/8,
\qquad 2(b+b_0)\sqrt n\,\varepsilon<1,
$$
and set
$$
\mathcal T=\{X\in\Sym(n):\dist(X,\Sigma_{n,k})\le\varepsilon\}.
$$
We claim that
\begin{equation}\label{eq:pucci-compatibility}
 \mathcal P^-_{a,b}(X-Y)
 \le H_0(X)-H_0(Y)
 \le\mathcal P^+_{a,b}(X-Y)
 \qquad\text{for all }X,Y\in\mathcal T.
\end{equation}
If \(|X-Y|<s/2\), every point on the segment from \(X\) to \(Y\)
is within distance \(\varepsilon+|X-Y|<s\) of \(\Sigma_{n,k}\).
Integrating \eqref{eq:local-operator-ellipticity} along this
segment proves \eqref{eq:pucci-compatibility}.

If \(|X-Y|\ge s/2\), choose nearest points
\(C,D\in\Sigma_{n,k}\) to \(X,Y\), respectively. Then
$$
|X-C|,\ |Y-D|\le\varepsilon,
\qquad |C-D|\ge s/2-2\varepsilon>s/4.
$$
Integrating along the segments from \(C\) to \(X\) and from
\(D\) to \(Y\), and using \(H_0(C)=H_0(D)=0\), we obtain
$$
|H_0(X)-H_0(Y)|\le2b_0\sqrt n\,\varepsilon.
$$
By \eqref{eq:far-pair-margin} and the choice of \(\varepsilon\),
$$
\mathcal P^+_{a,b}(X-Y)
\ge1-2b\sqrt n\,\varepsilon
>2b_0\sqrt n\,\varepsilon
\ge |H_0(X)-H_0(Y)|.
$$
The other inequality in \eqref{eq:pucci-compatibility} follows
by interchanging \(X\) and \(Y\) and using
\(\mathcal P^-_{a,b}(E)=-\mathcal P^+_{a,b}(-E)\).

The bounds in \eqref{eq:pucci-compatibility} allow us to use
the formula
$$
\widehat H(N)=\inf_{Y\in\mathcal T}
\left\{H_0(Y)+\mathcal P^+_{a,b}(N-Y)\right\},
\qquad N\in\Sym(n),
$$
as in \cite[Remark~1.5]{CabreCaffarelli2003}. For \(N\in\mathcal T\), \eqref{eq:pucci-compatibility} shows
that every term is at least \(H_0(N)\), with equality at \(Y=N\).
Thus \(\widehat H=H_0\) on \(\mathcal T\).
For every \(P\ge0\) and \(Y\in\mathcal T\),
$$
\mathcal P^+_{a,b}(N-Y)+a\tr P
\le\mathcal P^+_{a,b}(N+P-Y)
\le\mathcal P^+_{a,b}(N-Y)+b\tr P.
$$
Adding \(H_0(Y)\) and taking the infimum over \(Y\in\mathcal T\),
we obtain
\begin{equation}\label{eq:global-ellipticity}
 a\tr P\le\widehat H(N+P)-\widehat H(N)\le b\tr P
 \qquad\text{for all }N\in\Sym(n),\quad P\ge0.
\end{equation}
\end{proof}

The next lemma shows what happens when we add a small positive
quadratic term to a possibly degenerate dual potential. To be specific, the inverse matrices remain admissible, and the lemma controls the resulting error term in the equation.

\begin{lemma}\label{lem:regularize-cone}
There exist constants \(\varepsilon_0>0\) and \(C_0<\infty\),
depending only on \(n,k,K\), such that for every
\(C\in\Sigma_{n,k}\) and \(0<\varepsilon<\varepsilon_0\), we have
\begin{equation}\label{eq:regularized-physical}
 A_\varepsilon=(C+\varepsilon I)^{-1}-\kappa I\in\Gamma_k^n,
 \qquad
 \bigl|F_{n,k}(A_\varepsilon)-1\bigr|\le C_0\varepsilon.
\end{equation}
\end{lemma}

\begin{proof}
Lemma~\ref{lem:operator} allows us to choose a neighborhood
of \(\Sigma_{n,k}\) where \(T_{k-1}>0\), \(H_0\) is smooth,
and \(|H_0|<1/2\). We then choose \(\varepsilon_0>0\) small enough that \(C+tI\) stays in this neighborhood whenever \(C\in\Sigma_{n,k}\) and \(0\le t\le\varepsilon_0\).

First take \(C\) in the set in \eqref{eq:physical-set} before
the closure. For \(0\le t\le\varepsilon_0\), set
\(A(t)=(C+tI)^{-1}-\kappa I\). Then
$$
\sigma_{k-1}(A(t))
=\frac{T_{k-1}(C+tI)}{\det(C+tI)}>0,
\qquad
\frac{\sigma_k(A(t))}{\sigma_{k-1}(A(t))}
=1-H_0(C+tI)>\frac12.
$$
The Newton--Maclaurin inequalities imply
\(A(t)\in\Gamma_k^n\) throughout the path.

For general \(C\in\Sigma_{n,k}\), fix \(\varepsilon>0\)
and choose \(C_j\to C\) from the set in
\eqref{eq:physical-set} before taking the closure.
The argument above applies to
\((C_j+\varepsilon I)^{-1}-\kappa I\), which converges to
\(A_\varepsilon\). Hence \(A_\varepsilon\in\overline\Gamma_k^n\).
Since \(C+\varepsilon I\) lies in the chosen neighborhood of
\(\Sigma_{n,k}\),
$$
\sigma_{k-1}(A_\varepsilon)
=\frac{T_{k-1}(C+\varepsilon I)}{\det(C+\varepsilon I)}>0,
\qquad
\frac{\sigma_k(A_\varepsilon)}{\sigma_{k-1}(A_\varepsilon)}
=1-H_0(C+\varepsilon I)>\frac12.
$$
Thus \(\sigma_k(A_\varepsilon)>0\), and the Newton--Maclaurin
inequalities imply \(A_\varepsilon\in\Gamma_k^n\).
Finally,
$$
F_{n,k}(A_\varepsilon)-1=-H_0(C+\varepsilon I),
\qquad H_0(C)=0.
$$
The error estimate in \eqref{eq:regularized-physical} now follows
from the uniform bound for \(DH_0\) on this neighborhood.
\end{proof}

The last lemma shows how the condition on the Hessians passes
to a locally uniform limit.
\begin{lemma}\label{lem:physical-limit}
Let \(w_j\) be smooth functions on an open set
\(\Omega\subset\R^n\), with \(D^2w_j\in\Sigma_{n,k}\) everywhere.
If \(w_j\to W\) locally uniformly, then
\(W\in C^{1,1}_{\mathrm{loc}}(\Omega)\),
\(0\le D^2W\le I\) almost everywhere, and
$$
D^2W\in\Sigma_{n,k}\qquad\text{almost everywhere in }\Omega.
$$
\end{lemma}

\begin{proof}
Since the functions \(w_j\) and \(|y|^2/2-w_j\) are locally convex in \(\Omega\), we know \(W\) is also locally convex
and semiconcave, belongs to \(C^{1,1}_{\mathrm{loc}}(\Omega)\),
and satisfies \(0\le D^2W\le I\) almost everywhere.

At an Alexandrov point \(y_0\) of \(W\), set
$$
Q(y)=W(y_0)+DW(y_0)\cdot(y-y_0)
     +\frac12(y-y_0)^TA(y-y_0),
\qquad W(y)-Q(y)=o(|y-y_0|^2).
$$
For any fixed \(\varepsilon>0\), first choose \(r>0\) so that
\(\overline B_r(y_0)\subset\Omega\) and
$$
\sup_{\overline B_r(y_0)}|W-Q|\le\frac{\varepsilon r^2}{8}.
$$
Then choose \(j\) sufficiently large that
$$
\sup_{\overline B_r(y_0)}|w_j-W|
\le\frac{\varepsilon r^2}{8}.
$$
It follows that on \(\partial B_r(y_0)\),
$$
w_j-Q-\varepsilon|y-y_0|^2
\le-\frac{3\varepsilon r^2}{4},
$$
while its value at \(y_0\) is at least \(-\varepsilon r^2/8\).
Thus the maximum is attained at an interior point \(y_j^+\),
where
$$
C_j^+:=D^2w_j(y_j^+)\le A+2\varepsilon I.
$$
Similarly, \(w_j-Q+\varepsilon|y-y_0|^2\) is at least
\(3\varepsilon r^2/4\) on the boundary and at most
\(\varepsilon r^2/8\) at the center. At an interior minimum
\(y_j^-\),
$$
C_j^-:=D^2w_j(y_j^-)\ge A-2\varepsilon I.
$$

Take a sequence \(\varepsilon\downarrow0\), choosing the radius
and index for each \(\varepsilon\) as above. By compactness,
we can pass to a subsequence of the pairs
\((C_j^+,C_j^-)\). Their limits satisfy
$$
C^+\le A\le C^-,\qquad C^+,C^-\in\Sigma_{n,k}.
$$
The order property in Lemma~\ref{lem:operator} then implies
\(C^+=C^-\), so \(A=C^+=C^-\in\Sigma_{n,k}\).
Since Alexandrov points have full measure, the lemma follows.
\end{proof}

\section{The Jacobi inequality and a rescaled limit}
\label{sec:jacobi}
In this section, we study how the unboundedness of \(D^2u\) affects the dual
potential on expanding balls.
We first derive a Jacobi inequality for the logarithm of the largest
eigenvalue using the concavity inequality in \cite{DZ}.
\begin{lemma}\label{lem:jacobi}
There are constants $L_0>1$ and $c>0$, depending only on $n,k,K$, such that $b=\log\lambda_{\max}(D^2u)$ satisfies
\begin{equation}\label{eq:jacobi}
 F^{ij}\varphi_{ij}\ge cF^{ij}\varphi_i\varphi_j
\end{equation}
for every upper test $\varphi$ for $b$ at a point where $\lambda_{\max}\ge L_0$. Here $F^{ij}:=\frac{\partial F_{n,k}(D^2u)}{\partial u_{ij}}$.
\end{lemma}
\begin{proof}
Let $\varphi$ be an upper test for $b$ at $x_0$ and choose orthonormal
coordinates such that $D^2u(x_0)$ is diagonal. Suppose that
\[
 \lambda_1=\cdots=\lambda_m>\lambda_{m+1}\ge\cdots\ge\lambda_n.
\]
All derivatives below are evaluated at $x_0$. By \cite[Lemma 5]{BCD}, we have
\begin{equation}\label{eq:root-first-contact}
 u_{abi}=\delta_{ab}\lambda_1\varphi_i,
 \qquad a,b\le m,\quad 1\le i\le n,
\end{equation}
and
\begin{equation}\label{eq:root-second-contact}
 \varphi_i=\frac{u_{11i}}{\lambda_1},\qquad
 \varphi_{ii}\ge\frac{u_{11ii}}{\lambda_1}
 +\frac2{\lambda_1}\sum_{p>m}
       \frac{u_{1pi}^2}{\lambda_1-\lambda_p}
 -\frac{u_{11i}^2}{\lambda_1^2}.
\end{equation}
Setting $\zeta_p=u_{pp1}$ and differentiating the equation $F(D^2u)=1$ once and twice with respect to $x_1$, we obtain
\begin{equation}\label{eq:pde-third-1}
 \sum_pF_{p}\zeta_p=0,
\end{equation}
and
\begin{align}\label{eq:pde-third-2}
 \sum_iF_i u_{ii11}
 \ge-\sum_{p,q}F_{pq}\zeta_p\zeta_q
 +2\sum_{p>m}
       \frac{F_p-F_1}{\lambda_1-\lambda_p}u_{11p}^2.
\end{align}
Combining these with \eqref{eq:root-second-contact}, we get
\begin{align}
 \sum_iF_i\varphi_{ii}
 \ge\frac1{\lambda_1}
\left(-\sum_{p,q}F_{pq}\zeta_p\zeta_q
+2\sum_{p>m}
\frac{F_p\zeta_p^2}{\lambda_1-\lambda_p}\right)+\sum_{p>m}\frac{F_p}{\lambda_1^2}
\frac{\lambda_1+\lambda_p}{\lambda_1-\lambda_p}
u_{11p}^2-\frac{F_1u_{111}^2}{\lambda_1^2}.
\label{eq:jacobi-combination}
\end{align}
Let $L_0$ be sufficiently large so that
\[
L_0\geq \max\{L_*,K/\delta,3K\},
\]
where $L_*$, $\delta$, and $\gamma$ are as in \eqref{eq:tangent}.
Then, for $\lambda_1\geq L_0$, semiconvexity implies that
\begin{equation}\label{eq:large-root-gaps}
 \lambda_1-\lambda_p\le(1+\delta)\lambda_1,
 \qquad
 \frac{\lambda_1+\lambda_p}{\lambda_1-\lambda_p}\ge\frac12
 \qquad(p>m).
\end{equation}
Since $\zeta_p=0$ for $1<p\le m$, we may apply
\eqref{eq:tangent} to obtain
\[
 -\sum_{p,q}F_{pq}\zeta_p\zeta_q
 +2\sum_{p>m}\frac{F_p\zeta_p^2}{\lambda_1-\lambda_p}
 \ge(1+\gamma)\frac{F_1u_{111}^2}{\lambda_1}.
\]
Substituting into \eqref{eq:jacobi-combination} and using
\eqref{eq:large-root-gaps}, we conclude that
\[
 \sum_iF_i\varphi_{ii}
 \ge\gamma\frac{F_1u_{111}^2}{\lambda_1^2}
    +\frac12\sum_{p>m}\frac{F_pu_{11p}^2}{\lambda_1^2}
 \ge\min\{\gamma,\tfrac12\}\sum_iF_i\varphi_i^2.
\]
This proves \eqref{eq:jacobi} at every upper contact.
\end{proof}
Following Shankar and Yuan \cite{SY,SYrigidity}, we combine the
Jacobi inequality with the Legendre transform and the weak Harnack inequality. Inspired by Ding \cite{Ding}, we then use the unboundedness of $D^2u$ to show that the averages of $\det D^2w$ over expanding balls tend to zero. We use this decay to construct the limit needed for the splitting argument in the next section.
\begin{lemma}\label{lem:blowdown}
Let $v=u+\kappa|x|^2/2$ and $w=v^*$. If $\sup_{\R^n}\lambda_{\max}(D^2u)=\infty$, then there exists
a sequence \(R_j\to\infty\) such that the functions
\begin{equation}\label{eq:dual-rescaling}
 w_j(y)=\frac{w(R_jy)-w(0)-R_jDw(0)\cdot y}{R_j^2}
\end{equation}
converge locally in $C^1$ to a convex $W\in C^{1,1}_{\mathrm{loc}}(\R^n)$, with globally Lipschitz gradient, satisfying
\begin{equation}\label{eq:limit-properties}
 W(0)=0,\qquad DW(0)=0,\qquad
 D^2W\in\Sigma_{n,k},\qquad \det D^2W=0
 \quad\text{almost everywhere}.
\end{equation}
\end{lemma}
\begin{proof}
Following Shankar and Yuan \cite[Proposition~2.3]{SY}, define
\[
 G(N)=-F(N^{-1}-\kappa I)
 \qquad\text{for }N>0\text{ with }N^{-1}-\kappa I\in\Gamma_k^n.
\]
Under the global diffeomorphism $y=Dv(x)$ from
Section~\ref{sec:spectral}, the Legendre transform $w=v^*$ satisfies
\[
 D^2w(y)=(D^2u(x)+\kappa I)^{-1},
 \qquad G(D^2w)=-1.
\]
Differentiating $G$ with respect to its matrix argument, we obtain
\begin{equation}\label{eq:dual-linearization}
 G^{ij}(D^2w(y))
 =F^{pq}(D^2u(x))v_{pi}(x)v_{qj}(x),
 \qquad x=x(y).
\end{equation}
At a point where $D^2u$ is diagonal, the matrix
$(G^{ij}(D^2w))$ is diagonal with entries
\[
 G^{ii}(D^2w)=F_i(\lambda_i+\kappa)^2.
\]
Thus \eqref{eq:dual-ellipticity} ensures that this coefficient
matrix is uniformly elliptic throughout $\R^n$.
Below, $F^{pq}$ and $G^{ij}$ are evaluated at $D^2u(x)$ and
$D^2w(y)$, respectively.

For a smooth function $\phi$ of $y$, the transformation rule
in \cite[Proposition~2.3]{SY} implies
\[
 F^{pq}\partial_{x_p x_q}(\phi\circ Dv)
 =G^{ij}\phi_{ij}.
\]

Define
\[
 b(y)=\log\max\{\lambda_{\max}(D^2u(x(y))),2L_0\}.
\]
By Lemma~\ref{lem:jacobi} and the stability of viscosity
subsolutions under maxima, we have
\[
 G^{ij}b_{ij}\ge c\,G^{ij}b_i b_j
 \qquad\text{in }\R^n
\]
in the viscosity sense.

As in \cite[Proposition~2.2]{SYrigidity},
choose $0<\theta\le c$ and set
\begin{equation}\label{eq:jacobi-superfunction}
 q(y)=e^{-\theta b(y)}.
\end{equation}
It follows that \(q\) is positive and bounded, and it satisfies
\[
 G^{ij}q_{ij}\le0
 \qquad\text{in }\R^n
\]
in the viscosity sense.

The unboundedness assumption and the surjectivity of $y=Dv(x)$
imply $\inf_{\R^n}q=0$.
Following the weak Harnack argument in
\cite[proof of Proposition~3.1, Step~1]{SYrigidity},
we apply Lemma~\ref{lem:harnack} on $B_{2R}$ to obtain
\begin{equation}\label{eq:average-zero}
 \left(\fint_{B_R}q^s\right)^{1/s}
 \le C\inf_{B_R}q\longrightarrow0
 \qquad(R\to\infty),
\end{equation}
where $s>0$ and $C$ depend only on the dimension and the
uniform ellipticity constants.
Boundedness of $q$ and H\"older's inequality then imply
\begin{equation}\label{eq:harnack--t}
 \fint_{B_R}q^t\longrightarrow0
 \qquad(R\to\infty)
 \quad\text{for every }t>0.
\end{equation}

The eigenvalues of $D^2w$ lie in $(0,1]$, and its smallest eigenvalue is $(\lambda_{\max}+\kappa)^{-1}$. Thus
\begin{equation}\label{eq:det-control}
 0\le\det D^2w
 \le\frac1{\lambda_{\max}+\kappa}
 \le Cq^{1/\theta}.
\end{equation}
Consequently, for each fixed $r>0$,
 \eqref{eq:harnack--t} implies that
\begin{equation}\label{eq:scaled-mass-zero}
 0\le \int_{B_r}\det D^2w_j(y)\,dy
 =R_j^{-n}\int_{B_{rR_j}}\det D^2w(z)\,dz\le C R_j^{-n}\int_{B_{rR_j}}q(z)^{1/\theta}\,dz\longrightarrow0.
\end{equation}

For any sequence $R_j\to\infty$, the normalization in \eqref{eq:dual-rescaling} and $0<D^2w_j\le I$ imply
\[
 |Dw_j(y)|\le|y|,\qquad 0\le w_j(y)\le\tfrac12|y|^2.
\]
By the Arzel\`a--Ascoli theorem and a diagonal argument,
a subsequence converges locally in \(C^1\) to a convex function
\(W\) satisfying \(0\le D^2W\le I\) almost everywhere. Every $D^2w_j$ lies in $\Sigma_{n,k}$, so Lemma~\ref{lem:physical-limit} proves the third assertion in \eqref{eq:limit-properties}.

By the weak continuity of Alexandrov Monge--Amp\`ere measures
under local uniform convergence of convex functions, the measures
\(\det D^2w_j\,dy\) converge weakly to the Monge--Amp\`ere measure
of \(W\); see \cite{Gutierrez}.
It follows from \eqref{eq:scaled-mass-zero} that this limiting
measure is zero. Since \(DW\) is locally Lipschitz, the
Monge--Amp\`ere measure of \(W\) has density \(\det D^2W\),
and hence \(\det D^2W=0\) almost everywhere.
The normalization also passes to the limit under local \(C^1\)
convergence, so \(W(0)=0\) and \(DW(0)=0\).
\end{proof}

\section{Reduction to a lower-order equation}
\label{sec:splitting}
In this section, we construct a smooth admissible solution of a
lower-order equation from the limit obtained in
Section~\ref{sec:jacobi}. We first remove its global affine
directions and show that the conjugate of the remaining function
is finite everywhere. Regularization and interior estimates then
allow us to recover the required solution. Throughout this section,
the Legendre conjugate of a finite convex function $W$ is
\[
 W^*(p)=\sup_{y\in\R^n}\{p\cdot y-W(y)\},
 \qquad
 \dom W^*=\{p\in\R^n:W^*(p)<+\infty\}.
\]
We will use the following subdifferential identity
(see \cite{Rockafellar}):
\begin{equation}\label{eq:fenchel-reciprocity}
 p\in\partial W(y)
 \quad\Longleftrightarrow\quad
 y\in\partial W^*(p)
 \quad\Longleftrightarrow\quad
 W(y)+W^*(p)=p\cdot y.
\end{equation}

A nonzero vector $e$ is called a global affine direction of $W$ if
there is a constant $a$ such that
\[
 W(y+te)=W(y)+ta
 \qquad\text{for every }y\in\R^n\text{ and }t\in\R.
\]
The next lemma identifies the global affine directions of \(W\)
using the effective domain of its conjugate.

\begin{lemma}\label{lem:MA-split}
Let $W\in C^{1,1}_{\mathrm{loc}}(\R^n)$ be convex, with
\[
 W(0)=0,\qquad DW(0)=0,\qquad
 D^2W\in\Sigma_{n,k},\qquad \det D^2W=0
 \quad\text{almost everywhere}.
\]
After an orthogonal change of coordinates, there are an integer
$1\le m\le k-1$ and a convex function
$V\in C^{1,1}_{\mathrm{loc}}(\R^{n-m})$ such that
\begin{equation}\label{eq:maximal-splitting}
 W(y'',y')=V(y'),
 \qquad (y'',y')\in\R^m\times\R^{n-m}.
\end{equation}
The function $V$ has no nonzero global affine direction. Moreover,
\begin{equation}\label{eq:compressed-physical-spectrum}
 V(0)=0,\qquad DV(0)=0,\qquad
 0\le D^2V\le I,\qquad
 D^2V\in\Sigma_{n-m,k-m}
 \quad\text{almost everywhere}.
\end{equation}
\end{lemma}

\begin{proof}
Since \(W\in C^{1,1}_{\mathrm{loc}}(\R^n)\) and
\(\det D^2W=0\) almost everywhere, the area formula implies
$$
|DW(B_R)|=0
\qquad\text{for every }R>0.
$$
It follows that
\begin{equation}\label{eq:gradient-image-null}
 |DW(\R^n)|=0.
\end{equation}

Suppose that \(\dom W^*\) has full affine dimension.
By convexity, its interior is nonempty. For any \(p\) in this
interior, choose \(y\in\partial W^*(p)\).
Then \eqref{eq:fenchel-reciprocity} implies \(p=DW(y)\), so
$$
\operatorname{int}(\dom W^*)\subset DW(\R^n),
$$
contrary to \eqref{eq:gradient-image-null}. Thus the affine hull
of \(\dom W^*\) is a proper affine subspace of \(\R^n\).

Convexity and the normalization at the origin imply \(W\ge0\)
and \(W^*(0)=0\). Thus \(0\in\dom W^*\), and its affine hull
is the linear space
$$
L=\operatorname{span}(\dom W^*).
$$
Set \(E=L^\perp\) and \(m=\dim E\ge1\). Then
$$
W(y)=\sup_{p\in\dom W^*}\{p\cdot y-W^*(p)\}=W(P_Ly),
$$
where \(P_L\) denotes the orthogonal projection onto \(L\).
In particular, \(W(y+e)=W(y)\) for all \(e\in E\).
Defining \(V=W|_L\), we obtain \eqref{eq:maximal-splitting}.

We now show that \(V\) has no nonzero global affine direction.
Let \(e\in L\) and \(a\in\R\) satisfy
\(V(y'+te)=V(y')+ta\) for all \(y'\in L\) and \(t\in\R\).
Since \(DV(0)=0\), we have \(a=0\), and hence \(W(te)=0\)
for all \(t\in\R\). For every \(p\in\dom W^*\),
$$
\sup_{t\in\R}t(p\cdot e)\le W^*(p)<\infty,
$$
so \(p\cdot e=0\). Hence \(e\) is orthogonal to
\(L=\operatorname{span}(\dom W^*)\). Since \(e\in L\),
we conclude that \(e=0\).

At every point \(y'\) where \(V\) has a second-order expansion,
\eqref{eq:maximal-splitting} implies
\begin{equation}\label{eq:split-hessian-block}
 D^2W(y'',y')=0_m\oplus D^2V(y')
\end{equation}
for all \(y''\in E\). Since \(D^2W\in\Sigma_{n,k}\) almost
everywhere, Fubini's theorem implies
\(0_m\oplus D^2V(y')\in\Sigma_{n,k}\) for almost every \(y'\).
Applying Lemma~\ref{lem:strata} to this block matrix, we obtain
\(m\le k-1\) and \(D^2V\in\Sigma_{n-m,k-m}\) almost everywhere.
Together with the normalization and \(0\le D^2W\le I\),
this proves \eqref{eq:compressed-physical-spectrum}.
\end{proof}

The splitting \eqref{eq:maximal-splitting} implies
\begin{equation}\label{eq:split-conjugates}
 W^*(p'',p')=
 \begin{cases}
  V^*(p'),&p''=0,\\
  +\infty,&p''\ne0.
 \end{cases}
\end{equation}
Hence
$$
\dom W^*=\{0\}\times\dom V^*.
$$
To recover an entire solution in the remaining variables by
conjugation, we still need to prove
$$
\dom V^*=L.
$$

\begin{lemma}\label{lem:full-conjugate}
Let \(\mathcal H:\Sym(d)\to\R\) be continuous and
uniformly elliptic, and let
\(V\in C^{1,1}_{\mathrm{loc}}(\R^d)\) be convex.
Suppose that
\begin{equation*}
 \mathcal H(D^2V)=0
 \qquad\text{in }\R^d
\end{equation*}
in the viscosity sense. If \(V\) has no nonzero global affine
direction, then \(\dom V^*=\R^d\).
\end{lemma}

\begin{proof}
Suppose that \(D=\dom V^*\) is a proper subset of \(\R^d\).
By convexity and the separation
theorem, there exists a nonzero vector \(e\) such that
\begin{equation}\label{eq:bounded-effective-slopes}
 \sup_{p\in D}p\cdot e<\infty.
\end{equation}
Since \(DV(y)\in D\) for every \(y\in\R^d\), the function
$$
f(y)=V(y+e)-V(y)
    =\int_0^1 DV(y+te)\cdot e\,dt
$$
is bounded above on \(\R^d\).

Let \(\mathcal P^\pm\) be the Pucci operators associated with
the ellipticity constants of \(\mathcal H\).
Since \(V(\,\cdot+e)\) and \(V\) satisfy the same equation
almost everywhere, uniform ellipticity implies
\begin{equation}\label{eq:difference-pucci}
 \mathcal P^-(D^2f)\le0\le\mathcal P^+(D^2f)
 \qquad\text{almost everywhere in }\R^d.
\end{equation}
As \(f\in C^{1,1}_{\mathrm{loc}}(\R^d)\),
Lemma~\ref{lem:strong-viscosity}, applied to \(f\) and \(-f\),
shows that these inequalities also hold in the viscosity sense.

Set \(M=\sup_{\R^d}f<\infty\). Since \(M-f\) is nonnegative
and satisfies \eqref{eq:difference-pucci},
Lemma~\ref{lem:harnack} applied on \(B_{2R}\) yields
$$
\sup_{B_R}(M-f)
\le C\inf_{B_R}(M-f)
=C\left(M-\sup_{B_R}f\right)\longrightarrow0
\qquad(R\to\infty),
$$
where \(C\) is independent of \(R\). Hence \(f\equiv M\).

For each \(y\in\R^d\), the function
\(t\mapsto V(y+te)-Mt\) is convex and \(1\)-periodic,
hence constant. Thus
$$
V(y+te)=V(y)+Mt
\qquad(y\in\R^d,\ t\in\R).
$$
This contradicts the assumption that \(V\) has no nonzero
global affine direction. Therefore, we conclude that \(\dom V^*=\R^d\).
\end{proof}

Let \(V\) be the function from Lemma~\ref{lem:MA-split},
and set \(d=n-m\) and \(p=k-m\).
By Lemma~\ref{lem:operator}, there is a globally uniformly
elliptic operator \(\widehat H\) vanishing on \(\Sigma_{d,p}\).
Thus \eqref{eq:compressed-physical-spectrum} implies
\(\widehat H(D^2V)=0\) almost everywhere.
By Lemma~\ref{lem:strong-viscosity}, this equation also holds
in the viscosity sense. Lemma~\ref{lem:full-conjugate} then
implies that \(V^*\) is finite on all of \(\R^d\).

To recover a smooth admissible solution of the lower-order
equation from \(V^*\), we first prove the following regularity
result for admissible viscosity solutions.

\begin{lemma}\label{lem:visc-reg}
Every continuous admissible viscosity solution of
\[
 F_{d,p}(D^2z)=1
\]
in an open set $\Omega\subset\R^d$ is smooth in $\Omega$ and
satisfies $D^2z\in\Gamma_p^d$ there.
\end{lemma}

\begin{proof}
For \(p=1\), the conclusion follows from interior regularity
for the Poisson equation. Let \(p\ge2\) and fix a ball
\(\overline{B_R}\subset\Omega\).
Choose \(h_j\in C^\infty(\overline{B_R})\) converging uniformly
to \(z\) on \(\partial B_R\).
The Dirichlet existence theorem in \cite{CNS} yields smooth
admissible solutions \(z_j\) of
\begin{equation}\label{eq:smooth-dirichlet-approximants}
 F_{d,p}(D^2z_j)=1\quad\text{in }B_R,
 \qquad z_j=h_j\quad\text{on }\partial B_R.
\end{equation}

Set \(e_j=\|h_j-z\|_{L^\infty(\partial B_R)}\).
Comparing \(z\) with the strict perturbations
\(z_j\pm e_j\pm\delta(R^2-|x|^2)\) in the viscosity sense
and letting \(\delta\downarrow0\), we obtain
\begin{equation}\label{eq:dirichlet-uniform-convergence}
 \|z_j-z\|_{L^\infty(B_R)}
 \le e_j\longrightarrow0.
\end{equation}

The oscillations of \(z_j\) are uniformly bounded, so Chen's
interior gradient estimate
\cite[Theorem~1.1, equation~(1.3)]{Chen} yields a uniform
gradient bound on \(B_{3R/4}\).
The interior Hessian estimate of Qiu and Yan
\cite[Theorem~1.1]{QY} then yields
$$
\sup_{B_{R/2}}|D^2z_j|\le C,
$$
where \(C\) is independent of \(j\).

By Lemma~\ref{lem:cone-newton}, the Hessians \(D^2z_j\) lie
in a fixed compact subset of \(\Gamma_p^d\).
On a neighborhood of its convex hull, \(F_{d,p}\) is smooth,
concave, and uniformly elliptic.
The Evans--Krylov and Schauder estimates therefore yield
uniform bounds for all derivatives on smaller balls;
see \cite{CC}.
A subsequence converges smoothly on \(B_{R/4}\), and
\eqref{eq:dirichlet-uniform-convergence} identifies its limit
as \(z\). Since \(B_R\) was arbitrary, the conclusion holds throughout \(\Omega\).
\end{proof}

We next show that \(V^*-\kappa|x|^2/2\) satisfies the lower-order
equation. Since \(D^2V\) may be singular, we first regularize \(V\)
and then take its conjugate, using Lemma~\ref{lem:regularize-cone}
to control the error in the equation. Smoothness then follows from Lemma~\ref{lem:visc-reg}.
\begin{lemma}\label{lem:Moreau}
Let \(V\in C^{1,1}_{\mathrm{loc}}(\R^d)\) be convex and satisfy
$$
D^2V\in\Sigma_{d,p}\quad\text{almost everywhere},
\qquad \dom V^*=\R^d.
$$
Then the function
\begin{equation}\label{eq:recovered-entire-solution}
 z(x)=V^*(x)-\frac{\kappa}{2}|x|^2
\end{equation}
is a smooth entire admissible solution of
\(F_{d,p}(D^2z)=1\) and satisfies \(D^2z\ge-KI\).
\end{lemma}

\begin{proof}
By the definition of \(\Sigma_{d,p}\), we have
\(0\le D^2V\le I\) almost everywhere.
For \(\varepsilon>0\), set
\begin{equation}\label{eq:regularized-conjugates}
 V_\varepsilon(y)=V(y)+\frac{\varepsilon}{2}|y|^2,
 \qquad Z_\varepsilon=V_\varepsilon^*,
 \qquad z_\varepsilon(x)=Z_\varepsilon(x)-\frac{\kappa}{2}|x|^2.
\end{equation}
Since \(\varepsilon I\le D^2V_\varepsilon\le(1+\varepsilon)I\)
almost everywhere, \(DV_\varepsilon\) is a global bi-Lipschitz
homeomorphism with an \(\varepsilon^{-1}\)-Lipschitz inverse.
By \cite[Theorem~26.5]{Rockafellar},
\(DZ_\varepsilon=(DV_\varepsilon)^{-1}\), so
\(Z_\varepsilon\in C^{1,1}_{\mathrm{loc}}(\R^d)\).
Differentiating this identity yields almost everywhere
\begin{equation}\label{eq:ae-inverse-hessian}
 D^2Z_\varepsilon(x)
 =(D^2V(y)+\varepsilon I)^{-1}
 \ge\frac{1}{1+\varepsilon}I,
 \qquad x=DV_\varepsilon(y).
\end{equation}
The set of points $y$ where $DV$ is not differentiable or
$D^2V(y)\notin\Sigma_{d,p}$ has measure zero.
Its image under the Lipschitz map $DV_\varepsilon$ also has
measure zero, so for almost every $x$, the Hessian identity
holds with $D^2V(y)\in\Sigma_{d,p}$.

Lemma~\ref{lem:regularize-cone} therefore implies, for all
sufficiently small $\varepsilon$,
\begin{equation}\label{eq:ae-approximate-lower-equation}
 D^2z_\varepsilon\in\Gamma_p^d,
 \qquad
 \bigl|F_{d,p}(D^2z_\varepsilon)-1\bigr|\le C\varepsilon
 \quad\text{almost everywhere},
\end{equation}
with $C$ independent of $\varepsilon$ and of the point.
Choose $\varepsilon$ small enough that $1-C\varepsilon>0$, and
consider the sets
\[
 \mathcal E_\varepsilon
 =\{A\in\Gamma_p^d:F_{d,p}(A)\ge1-C\varepsilon\},
\]
and
\[
 \mathcal L_\varepsilon
 =\bigl(\Sym(d)\setminus\Gamma_p^d\bigr)
 \cup\{A\in\Gamma_p^d:F_{d,p}(A)\le1+C\varepsilon\}.
\]
The sets $\mathcal E_\varepsilon$ and $\mathcal L_\varepsilon$ are closed in $\Sym(d)$. Moreover, ellipticity and the monotonicity of $\Gamma_p^d$ imply
\[
 A\in\mathcal E_\varepsilon
 \ \Longrightarrow\ A+P\in\mathcal E_\varepsilon,
 \qquad
 A\in\mathcal L_\varepsilon
 \ \Longrightarrow\ A-P\in\mathcal L_\varepsilon
 \qquad(P\ge0).
\]

Since $z_\varepsilon\in C^{1,1}_{\mathrm{loc}}$, both
$z_\varepsilon$ and $-z_\varepsilon$ are locally semiconvex.
By \eqref{eq:ae-approximate-lower-equation}, their Hessians
belong almost everywhere to $\mathcal E_\varepsilon$ and
$-\mathcal L_\varepsilon$, respectively. Lemma \ref{lem:strong-viscosity}
therefore implies that
\[
 1-C\varepsilon\le F_{d,p}(D^2z_\varepsilon)
 \le1+C\varepsilon
\]
in the admissible viscosity sense. 

By definition,
\begin{equation}\label{eq:conjugate-monotone-limit}
 Z_\varepsilon(x)
 =\sup_y\left\{x\cdot y-V(y)-\frac\varepsilon2|y|^2\right\}
 \uparrow V^*(x)
 \qquad\text{as }\varepsilon\downarrow0.
\end{equation}
Since $V^*$ is finite and convex on $\R^d$, it is continuous.
Dini's theorem therefore implies $z_\varepsilon\to z$
locally uniformly. Stability of viscosity inequalities yields
\begin{equation}\label{eq:limiting-viscosity-equation}
 F_{d,p}(D^2z)=1
 \qquad\text{in }\R^d
\end{equation}
in the admissible viscosity sense. 

Finally, \eqref{eq:ae-inverse-hessian} shows that
$Z_\varepsilon-|x|^2/(2(1+\varepsilon))$ is convex.
Passing to the locally uniform limit and using $\kappa=K+1$,
we conclude that
\[
 V^*(x)-\frac12|x|^2
 =z(x)+\frac K2|x|^2
\]
is convex. Lemma~\ref{lem:visc-reg} now shows that $z$ is
smooth and admissible, and the preceding convexity implies
$D^2z\ge-KI$ pointwise.
\end{proof}

\section{Proof of the main theorem}\label{sec:main-proof}
In this section, we prove the rigidity results. If \(D^2u\) is unbounded, the preceding results produce a smooth
entire solution of lower order with the same semiconvexity bound.
The induction hypothesis makes the dual limit quadratic, and
Savin's theorem then yields a contradiction. We first settle the bounded-Hessian case by a scaling argument.
\begin{lemma}\label{lem:bounded-rigidity}
Let $u$ be a smooth entire admissible solution of $F_{n,k}(D^2u)=1$. If $D^2u$ is bounded, then $u$ is quadratic.
\end{lemma}
\begin{proof}
For $k=1$, every $u_{ij}$ is a bounded entire harmonic function, hence constant. We may therefore assume that $k\ge2$.

For $R\ge1$, define
\[
 u_R(x)=R^{-2}\bigl(u(Rx)-u(0)-RDu(0)\cdot x\bigr).
\]
It solves the same equation, has the same Hessian bound, and satisfies $\|u_R\|_{L^\infty(B_1)}\le C$ by Taylor expansions. The Evans--Krylov estimate therefore implies
\[
 [D^2u_R]_{C^\alpha(B_{1/2})}\le C
\]
for fixed $\alpha>0$ and $C$, independent of $R$. Equivalently, for fixed $x,y\in\R^n$ and $R>2\max\{|x|,|y|\}$, we have
\[
 |D^2u(x)-D^2u(y)|\le C\left(\frac{|x-y|}{R}\right)^\alpha.
\]
Letting \(R\to\infty\), we conclude that \(D^2u\) is constant. This completes the proof.
\end{proof}

We are now ready to prove the main theorem.
\begin{proof}[Proof of Theorem~\ref{thm:main}]
As in \cite{QY}, we reduce the case \(k-l=2\) to an adjacent
quotient equation by setting
$$
U(x,t)=u(x)+\frac12t^2.
$$
Let \(\lambda=\lambda(D^2u)\). Then the eigenvalues of \(D^2U\)
are \((\lambda,1)\). For \(1\le j\le k\), we have
$$
\sigma_j(\lambda,1)
=\sigma_j(\lambda)+\sigma_{j-1}(\lambda)>0.
$$
The equation for \(u\) also implies
$$
\sigma_k(\lambda,1)-\sigma_{k-1}(\lambda,1)
=\sigma_k(\lambda)-\sigma_{k-2}(\lambda)=0.
$$
Thus \(U\) is an entire \(k\)-admissible solution of
\(F_{n+1,k}(D^2U)=1\), with \(D^2U\ge-KI\).
Rigidity for adjacent quotients therefore
implies that \(U\), and hence \(u\), is quadratic.

We next prove rigidity for the adjacent quotient equation by
induction on $k$. For $k=1$, the equation reduces to
\[
\Delta u=\sum_{i=1}^n\lambda_i=1,
\]
where $\lambda_1,\ldots,\lambda_n$ are the eigenvalues of $D^2u$.
By the semiconvexity assumption $D^2u\geq-KI$, we have
$\lambda_i\geq-K$ for every $i=1,\ldots,n$. Consequently,
\[
\lambda_i
=1-\sum_{j\neq i}\lambda_j
\leq 1+(n-1)K.
\]
Thus,
\[
-K\leq\lambda_i\leq 1+(n-1)K,
\qquad i=1,\ldots,n,
\]
which implies that $D^2u$ is uniformly bounded in $\mathbb R^n$.
By Lemma~\ref{lem:bounded-rigidity}, $u$ is a quadratic polynomial.

Let \(2\le k\le n\), and assume that the result holds for every
\(1\le p<k\). By Lemma~\ref{lem:bounded-rigidity}, it remains to show that \(\lambda_{\max}(D^2u)\) is bounded.
If it were unbounded, Lemma~\ref{lem:blowdown} would yield
normalized rescalings \(w_j\) converging locally in \(C^1\)
to a function \(W\) satisfying \eqref{eq:limit-properties}.

Lemma~\ref{lem:MA-split} removes all fixed affine directions. After an orthogonal change of coordinates,
\begin{equation}\label{eq:induction-split}
 W(y'',y')=V(y'),\qquad
 (y'',y')\in\R^m\times\R^{n-m},\qquad 1\le m\le k-1,
\end{equation}
where $V$ has no nonzero global affine direction and 
\begin{equation}\label{eq:D2VinSigma}
 D^2V\in\Sigma_{n-m,k-m}\quad\hbox{almost everywhere}.
\end{equation}
Write $d=n-m$ and $p=k-m$. By Lemma~\ref{lem:operator} and
\eqref{eq:D2VinSigma}, $V$ satisfies
\[
 \widehat H_{d,p}(D^2V)=0
 \qquad\text{almost everywhere in }\R^d,
\]
where $\widehat H_{d,p}$ is globally uniformly elliptic. By Lemma~\ref{lem:strong-viscosity}, $\widehat H_{d,p}(D^2V)=0$ also in the viscosity sense. Moreover, Lemma~\ref{lem:full-conjugate} implies $\dom V^*=\R^d$.

By Lemma~\ref{lem:Moreau}, the function
\[
 z(x)=V^*(x)-\frac\kappa2|x|^2
\]
is a smooth entire admissible solution of
$F_{d,p}(D^2z)=1$ with $D^2z\ge-KI$.
Since $1\le p<k$, the induction hypothesis implies that
$z$ is quadratic. Thus $V^*=z+\kappa|x|^2/2$ is quadratic
with $D^2V^*\ge I$, and its conjugate $V$ is also quadratic.
Combining \eqref{eq:induction-split} with $W(0)=0$, $DW(0)=0$, and \eqref{eq:limit-properties}, we obtain
\begin{equation}\label{eq:quadratic-limit}
 W(y)=\frac12y^TAy,
 \qquad A\in\Sigma_{n,k},
 \qquad \det A=0.
\end{equation}

Finally, we show that $D^2w(0)=A$.
Let $\widehat H$ be the operator from
Lemma~\ref{lem:operator}, and define
\[
 G(N)=\widehat H(A+N),
 \qquad e_j=w_j-W.
\]
Since $A\in\Sigma_{n,k}$, the operator $G$ is globally
uniformly elliptic, is smooth near zero, and satisfies
$G(0)=0$. Its first and second derivatives are bounded
on a sufficiently small neighborhood of zero.
Moreover, $D^2W=A$ and the equation for $w_j$ imply
\[
 G(D^2e_j)=\widehat H(D^2w_j)=0,
 \qquad
 \|e_j\|_{L^\infty(B_1)}\longrightarrow0.
\]
Savin's small perturbation theorem
\cite[Theorem~1.3]{Savin} therefore yields
\[
 \|e_j\|_{C^{2,\alpha}(B_{1/2})}\le C
\]
for some $\alpha\in(0,1)$ and all sufficiently large $j$,
with $C$ independent of $j$.

By compactness, a subsequence of \(\{e_j\}\) converges to zero
in \(C^2(\overline{B_{1/4}})\). In particular,
\[
 D^2e_j(0)\longrightarrow0
\]
along this subsequence.
On the other hand, differentiating
\eqref{eq:dual-rescaling} twice shows that
$D^2w_j(0)=D^2w(0)$. Hence
\[
 D^2e_j(0)=D^2w(0)-A,
\]
and therefore $D^2w(0)=A$. This is impossible because
\begin{equation*}
 D^2w(0)
 =\bigl(D^2u(x(0))+\kappa I\bigr)^{-1}>0,
\end{equation*}
whereas \(\det A=0\).
Thus \(D^2u\) is bounded, and
Lemma~\ref{lem:bounded-rigidity} implies that \(u\) is quadratic.

This completes the induction and hence the proof.
\end{proof}

\section*{Acknowledgments}

The first and second authors acknowledge support from Grants~2025YFA1017603
(National Key R\&D Program of China) and~12571227
(National Natural Science Foundation of China).

\noindent The third and fourth authors acknowledge support from Grants 12661038
(National Natural Science Foundation of China).

\noindent The fourth author acknowledge support from Grants 12141103
(National Natural Science Foundation of China).

\section*{Declaration on the use of artificial intelligence}

The authors proposed using interior regularity estimates for
Hessian quotient equations, together with a Jacobi inequality
for the largest Hessian eigenvalue, to prove rigidity.
ChatGPT developed the remaining arguments and carried out
the detailed calculations.
The authors have independently checked and reorganized the
entire manuscript and take full responsibility for its contents.

\begingroup
\small
\raggedright
\bibliographystyle{amsplain}
\bibliography{references}
\endgroup
\printauthorinformation

\end{document}